\documentclass[12pt,amsmath,amscd,amsbsy,amssymb]{amsart}

\usepackage[cp1251]{inputenc}
\usepackage[english]{babel}
\usepackage{amssymb,amsmath,amsthm,amsxtra,amsfonts}

 \theoremstyle{plain}
\newtheorem{thm}{Theorem}[section]

\newtheorem{theorem}[thm]{Theorem}
\newtheorem{corollary}[thm]{Corollary}
\newtheorem{lemma}[thm]{Lemma}

\newenvironment{proposition*}[1]{\smallskip\noindent{\bf #1.}\it}{\medskip}
\newenvironment{theorem*}[1]{\smallskip\noindent{\bf #1.}\it}{\medskip}

\numberwithin{equation}{section} 
\renewcommand\Im{\operatorname{Im}}
\renewcommand\Re{\operatorname{Re}}
\newcommand{\Tr}{\operatorname{Tr}}

\newcommand{\R}{\mathbb{R}}

\title{$L^1$--spectrum of Banach space valued Ornstein--Uhlenbeck operators}
\author{Rostyslav V. Kozhan}
\begin{document}
\maketitle
\begin{abstract}
We characterize the $L^1(E,\mu_\infty)$--spectrum of the
Ornstein--Uhlenbeck ope\-rator $ Lf(x)=\frac12 \Tr{QD^2f(x)}+\langle
Ax,Df(x)\rangle $, where $\mu_\infty$ is the invariant measure for
the Ornstein--Uh\-len\-beck semigroup generated by $L$. The main
result covers the general case of an infinite-dimensional Banach
space $E$ under the assumption that the point spectrum of $A^*$ is
nonempty and extends several recent related results.\end{abstract}



\begin{section}{Introduction and results}\label{section1}

 In this paper we investigate spectral properties
of the generator of the transition semigroup associated with the
stochastic linear Cauchy problem
\begin{equation}\label{eq1}
\left\{
\begin{array}{rcl}
dU_t&=&AU_tdt+BdW^H_t,\\
U_0&=&x,
\end{array}
\right.
\end{equation}
where $A$ is the generator of a $C_0$-semigroup $(S(t))_{t\ge0}$ on
a real Banach space $E$, $B$ is a nonzero bounded operator from a
real Hilbert space $H$ into $E$, $(W_t^H)_{t\ge0}$ is an
$H$-cylindrical Wiener process, and $x\in E$. The problem is a
natural infinite-dimensional generalization of the Langevin equation
and arises in many applications, for example in optimal control
theory and interest rate models, see~\cite{{Pr,Z2}, {Pr,Z}}.

It is well known~\cite{{B,N}, {Pr,Z}} that the problem~(\ref{eq1})
admits a unique weak solution $\mathbf{U}=(U_t(x))_{t\ge0}$ if and
only if for all $t\in(0,\infty)$ there exists a centered Gaussian
Radon measure $\mu_t$ on E with covariance operator $Q_t\in
\mathcal{L}(E^*,E)$ given by
\begin{equation}\label{eq2}
\langle Q_tx^*,y^*\rangle=\int_0^t\langle
S(s)BB^*S^*(s)x^*,y^*\rangle\, ds,\qquad x^*,y^*\in E^*,
\end{equation}
and in this case the solution can be represented in the form
\begin{equation}\label{proc}
U_t(x)=S(t)x+\int_0^tS(t-s)B\,dW_t^H.
\end{equation}
The process $\mathbf{U}$ in~(\ref{proc}) is Gaussian and Markov and
its transition semigroup $\mathbf{P}=(P(t))_{t\ge0}$ (the
Ornstein--Uhlenbeck semigroup) is defined on bounded Borel functions
on $E$ by
\begin{equation*}
\left(P(t)f\right)(x):=\mathbb{E}\left(f(U_t(x))\right)=\int_E{f(S(t)x+y)\,d\mu_t(y)}.
\end{equation*}

Assume that the limit $Q_\infty:=\lim_{t\to\infty}Q_t$ exists in the
weak operator topology of $\mathcal{L}(E^*,E)$ and that there exists
a centered Gaussian Radon measure $\mu_\infty$ with covariance
operator $Q_\infty$. Under this assumption the measure $\mu_\infty$
is invariant for $\mathbf{P}$ (see~\cite{{B,N}, {Pr,Z}}), i.e.
\begin{equation*}
\int_EP(t)f(x)\,d\mu_\infty(x)=\int_Ef(x)\,d\mu_\infty(x),\qquad
t\ge0.
\end{equation*}
Throughout the paper we assume $\mu_\infty$ to be nondegenerate.

The above equality easily implies (see~\cite[Thm XIII.1]{Yos}) that
the semigroup $\mathbf{P}$ has a unique extension to a strongly
continuous contraction semigroup on $L^p(E,\mu_\infty)$,
$p\in[1,\infty)$, which we also denote by $\mathbf{P}$. Properties
of $\mathbf{P}$ in $L^p(E,\mu_\infty)$ for $p\in(1,\infty)$ have
been extensively investigated in the literature
(see~\cite{{Ch-M,G},{Ch-M,G2},{Pr,Z2},{M,P,Pr},N} and references
therein). Properties of $\mathbf{P}$ in $L^1(E,\mu_\infty)$ turn out
to be completely different from those in the spaces
$L^p(E,\mu_\infty)$, $p\in(1,\infty)$.
In particular for $p=1$ the semigroup $\mathbf{P}$ loses its
regularity properties, which it possesses in the case
$p\in(1,\infty)$ (see~\cite{{Ch-M,G2},{Pr,Z2},{G,N}}), and the
spectrum of its generator, which is $p$-independent for
$p\in(1,\infty)$ (see~\cite{{M,P,Pr},N}), changes drastically.


The key issue investigated in the present paper is the structure of
the spectrum of the generator $L$, called the Ornstein--Uhlenbeck
operator, of this semigroup in $L^1(E,\mu_\infty)$. Denote by
$\mathcal{C}(L)$ the set of continuous on $E$ functions of the form
$f(x)=\phi(\langle x,x_1^*\rangle,\ldots,\langle x,x_n^*\rangle)$,
where $x_j^*\in \mathcal{D}(A^*)$ for all $j=1,\ldots,n$ and
$\phi\in C^2(\mathbb{R}^n)$ with compact support. It follows
from~\cite{{G,N}} that $\mathcal{C}(L)$ is a core for $L$, and for
$f\in \mathcal{C}(L)$,
\begin{equation}\label{eq3}
Lf(x)=\frac12 \Tr{D_H^2f(x)}+\langle Ax,Df(x)\rangle,
\end{equation}
where $Df:E\to E^*$ is the Fr\'{e}chet derivative of $f$, and
$D_Hf:E\to H^*$ is the Fr\'{e}chet derivative of $f$ in the
direction of $H$ defined by
$$
D_Hf(x)=\sum_{j=1}^n \frac{\partial \phi}{\partial x_j}(\langle
x,x_1^*\rangle,\ldots,\langle x,x_n^*\rangle) B^*x_j^*.
$$
Denote $Q=BB^*\in \mathcal{L}(E^*,E)$. When $E$ is a Hilbert space
itself, the first term on the right-hand side of~(\ref{eq3}) becomes
just $\frac12 \Tr QD^2f(x)$.

By the result in~\cite[Thm 5.1]{{M,P,Pr}}, if $E=\mathbb{R}^n$ the
$L^1(\mathbb{R}^n,\mu_\infty)$--spectrum of $L$ is equal to
$\overline{\mathbb{C}}_-=\left\{\lambda: \Re\lambda\le0 \right\}$
with each $\lambda\in\mathbb{C}_-=\left\{\lambda: \Re\lambda<0
\right\}$ being an eigenvalue. This result was extended to
infinite-dimensional Banach spaces $E$ in~\cite{{vN,Pr}} under the
assumption of eventual compactness of the semigroup
$(S(t))_{t\ge0}$, and in~\cite{{Ch-M}} under the assumption that the
part of $A^*$ in the reproducing kernel Hilbert space of
$\mu_\infty$ has an eigenvalue $\gamma\in\mathbb{C}_-$.
Theorem~\ref{thm} of the present paper generalizes both of these
results while requiring less effort to prove it. Also it can serve
as an alternative simple coordinate-free proof of the corresponding
finite-dimensional result of~\cite{{M,P,Pr}}.

\begin{theorem}\label{thm}
If the point spectrum of $A^*$ is not empty
$\sigma_p(A^*)\ne\varnothing$, then the spectrum of the
Ornstein--Uhlenbeck operator $L$ coincides with
$\overline{\mathbb{C}}_-$, and each $\lambda\in\mathbb{C}_-$ is its
eigenvalue.
\end{theorem}


We remark that due to~\cite[Pr 2.5]{N-Un} we have
$\sigma_p(A^*)\cap\{\lambda: \Re{\lambda}\ge0\}=\varnothing$, and
thus the condition $\sigma_p(A^*)\ne\varnothing$ of
Theorem~\ref{thm} is equivalent to the existence of an eigenvalue of
$A^*$ with negative real part. An extension of Theorem~\ref{thm} for
the case $\sigma_p(A^*)=\varnothing$ seems to be an open question so
far.

\begin{corollary}
Under the assumptions of Theorem~\ref{thm}, the Ornstein-Uhlenbeck
semigroup $(P(t))_{t\ge0}$ on $L^1(E,\mu_\infty)$ is norm
discontinuous everywhere.

\begin{proof} Indeed, it easily follows from~\cite[Thm
4.18]{{En,N}} that the spectrum of the ge\-ne\-ra\-tor of an
eventually norm continuous semigroup cannot be equal to
$\overline{\mathbb{C}}_-$.\end{proof}
\end{corollary}
\end{section}

\emph{Acknowledgement.} The author would like to thank
J.M.A.M.~van~Neerven and Rostyslav~O.~Hryniv for helpful remarks.

\begin{section}{Proof of Theorem~\ref{thm}}\label{section2}

That each $\lambda\in\mathbb{C}_-$ is an eigenvalue of $L$ we
establish in Lemmas~\ref{lm1} and~\ref{lm2}. This implies
$\mathbb{C}_-\subset\sigma_p(L)\subset\sigma(L)$. The fact that
$\mathbf{P}$ is contractive on $L^1(E,\mu_\infty)$ implies
$\sigma(L)\subset\overline{\mathbb{C}}_-$. Since the spectrum is
closed, this finishes the proof.

By the argument in Section~\ref{section1}, we may assume that $A^*$
has an eigenvalue $\gamma\in\mathbb{C}_-$. Denote the corresponding
eigenvector as $x_0^*\in E^*_{\mathbb{C}}$, where $E^*_{\mathbb{C}}$
is the complexification of $E^*$.

\begin{lemma}\label{lm1}
If $\gamma\in\R\cap\mathbb{C}_-$, then each $\lambda\in\mathbb{C}_-$
is an eigenvalue of $L$.
\begin{proof}
Since $\gamma\in\R$, the corresponding eigenvector $x_0^*$ of $A^*$
may be chosen in $E^*$. We will show that for each
$\lambda\in\mathbb{C}_-$ we can find an eigenfunction of $L$ of the
form $f_\lambda(x)=\phi_\lambda(\langle x,x_0^*\rangle):E\rightarrow
\mathbb{C}$ with some function $\phi_\lambda$ on $\R$.

Consider the one-dimensional Ornstein--Uhlenbeck operator defined by
$$
L_1\phi(t)=\frac12 q\phi''(t)+\gamma t\phi'(t)
$$
for $\phi\in \mathcal{C}(L_1)=\{\phi\in C^2(\mathbb{R}) \mbox{ with
compact support}\}$, with $q=\langle Qx_0^*,x_0^*\rangle\ge0$. In
fact $\langle Qx_0^*,x_0^*\rangle\ne0$: $x_0^*$ is an eigenvector of
$S^*$, so (\ref{eq2}) and $\langle Qx_0^*,x_0^*\rangle=0$ imply
$\langle Q_\infty x_0^*,x_0^*\rangle=0$ which contradicts the
assumption that $\mu_\infty$ is nondegenerate. Hence $q>0$.

Here $\mathcal{C}(L_1)$ is viewed as a subspace of
$L^1(\mathbb{R},\nu_\infty)$, where $\nu_\infty$ is the invariant
measure for $L_1$.  Now observe that $\phi(t)\in L^1(\R,\nu_\infty)$
is equivalent to $\phi(\langle x,x_0^*\rangle)\in
L^1(E,\mu_\infty)$. Indeed, $\nu_\infty$ is a centered
one-dimensional Gaussian measure with variance
$\int_0^\infty{e^{\gamma s}q e^{\gamma
s}}ds=-\frac{1}{2\gamma}\langle Qx_0^*,x_0^*\rangle$. By definition,
measure $\mu_\infty$ on the cylindrical function $\phi(\langle
x,x_0^*\rangle)$ is centered one-dimensional Gaussian with variance
\begin{equation*}
\begin{split}
\langle Q_\infty x_0^*,x_0^*\rangle=\int_0^\infty{\langle
S(s)QS^*(s)x_0^*,x_0^*\rangle}ds=\int_0^\infty{\langle Qe^{\gamma
s}x_0^*,e^{\gamma s}x_0^*\rangle}ds=-\frac{1}{2\gamma}\langle
Qx_0^*,x_0^*\rangle.
\end{split}
\end{equation*}
Hence $\phi(t)\in L^1(\R,\nu_\infty)$ if and only if $\phi(\langle
x,x_0^*\rangle)\in L^1(E,\mu_\infty)$, and
$L^1(\R,\nu_\infty)$-convergence is equivalent to
$L^1(E,\mu_\infty)$-convergence for the corresponding $\langle
x,x_0^*\rangle$-cylindrical functions.

Now, for $\phi\in \mathcal{C}(L_1)$, we have $f(x):=\phi(\langle
x,x_0^*\rangle)\in \mathcal{C}(L)$ and
\begin{equation*}
\langle Ax,Df(x)\rangle=\langle Ax,\phi'(\langle
x,x_0^*\rangle)\,x_0^*\rangle=\phi'(\langle x,x_0^*\rangle)\,\langle
x,A^*x_0^*\rangle=\gamma\, \phi'(\langle x,x_0^*\rangle)\,\langle
x,x_0^*\rangle,
\end{equation*}
\begin{equation*}
\left(D_H^2f(x)\right)(y)=\phi''(\langle x,x_0^*\rangle)\, \langle
y,B^*x_0^*\rangle\,B^*x_0^*,\hspace{0.4cm} y\in H.
\end{equation*}
The only nonzero eigenvalue of the operator on the right-hand side
of the last equality is $\phi''(\langle x,x_0^*\rangle)\langle
Qx_0^*,x_0^*\rangle\, $, so
\begin{equation}\label{eq48}
Lf(x)=\frac12\langle Qx_0^*,x_0^*\rangle\, \phi''(\langle
x,x_0^*\rangle)+\gamma\,\langle x,x_0^*\rangle\, \phi'(\langle
x,x_0^*\rangle)=\left(L_1\phi\right)(\langle x,x_0^*\rangle).
\end{equation}

For each $\lambda\in\mathbb{C}_-$ let $\phi_\lambda(t)\in
L^1(\R,\nu_\infty)$ be an eigenfunction of $L_1$ corresponding to
the eigenvalue $\lambda$ (we use the one-dimensional case
of~\cite[Thm 5.1]{{M,P,Pr}}). Now we find $\phi_n(t)\in
\mathcal{C}(L_1)$ with $\phi_n\to \phi_\lambda$ in
$L^1(\R,\nu_\infty)$. Then $\phi_n(\langle x,x_0^*\rangle)\in
\mathcal{C}(L)$ and $\phi_n(\langle x,x_0^*\rangle)\to
\phi_\lambda(\langle x,x_0^*\rangle)$ in $L^1(E,\mu_\infty)$, and
applying~(\ref{eq48}) to $\phi_n(\langle x,x_0^*\rangle)$'s we
obtain $f_\lambda(x):=\phi_\lambda(\langle x,x_0^*\rangle)\in
\mathcal{D}(L)$ with $Lf_\lambda(x)=\lambda f_\lambda(x)$.

\end{proof}
\end{lemma}

We state an easy auxiliary lemma, the proof of which is left to the
reader.

\begin{lemma}\label{trace}
Let $H$ be a Hilbert space, $x_1,x_2,y_1,y_2\in H$. The trace of the
operator $Ax=\langle x,y_1\rangle x_1 +\langle x,y_2\rangle x_2$ is
equal to $\Tr{A}=\langle x_1,y_1\rangle +\langle x_2,y_2\rangle$.
\end{lemma}

Now we prove

\begin{lemma}\label{lm2}

Let $\gamma\in\mathbb{C}_-\setminus\R$. Then each
$\lambda\in\mathbb{C}_-$ is an eigenvalue of $L$.
\begin{proof}

Let $\gamma=a+bi$, $a<0$, $b\ne0$. Take $h_1^*:=\Re{x_0^*}\in E^*$,
$h_2^*:=\Im{x_0^*}\in E^*$. We have $A^*h_1^*=ah_1^*-bh_2^*$,
$A^*h_2=bh_1^*+ah_2^*$, and also
\begin{eqnarray}
S^*(s)h_1^*=e^{as}(h_1^*\cos{bs}-h_2^*\sin{bs}),\label{eq111}\\
S^*(s)h_2^*=e^{as}(h_1^*\sin{bs}+h_2^*\cos{bs}).\label{eq222}
\end{eqnarray}

We follow the same approach as in Lemma~\ref{lm1}. Consider the
two-dimensional Ornstein--Uhlenbeck operator
$$L_2\phi(t):=\frac12\Tr(RD^2\phi(t))+\langle Ct,D\phi(t)\rangle, \qquad t\in\R^2$$
for $\phi\in\mathcal{C}(L_2)=\{\phi\in C^2(\mathbb{R}^2) \mbox{ with
compact support}\}$, where $$R:=(r_{ij})_{i,j=1}^2=(\langle
Qh_i^*,h_j^*\rangle)_{i,j=1}^2,\quad
C:=(c_{ij})_{i,j=1}^2=\left(\begin{array}{cc}a&-b\\b&a\end{array}\right).$$
Now we show that the covariance operator $R_\infty$ of the invariant
measure $\nu_\infty$ corresponding to $L_2$ is the same as of the
image measure of $\mu_\infty$ under the map $(\langle
x,h_1^*\rangle, \langle x,h_2^*\rangle)$. This can be verified
directly: using~(\ref{eq111}) and (\ref{eq222}), it is easy to show
that the matrix $e^{Cs}Re^{C^*s}$ equals to the matrix
$(\langle S(s)QS^*(s)h_i^*,h_j^*\rangle)_{i,j=1}^2$, and thus
\begin{equation}\label{eq5}
\begin{split}
R_\infty=\int_0^{\infty}{e^{Cs}Re^{C^*s}}\,ds&=
\left(\int_0^{\infty}{\langle
S(s)QS^*(s)h_i^*,h_j^*\rangle}\,ds\right)_{i,j=1}^2\\&=
\left(\langle Q_\infty h_i^*,h_j^*\rangle \right)_{i,j=1}^2.
\end{split}
\end{equation}
This implies that $\phi(t)\in L^1(\R^2,\nu_\infty)$ if and only if
$\phi(\langle x,h_1^*\rangle, \langle x,h_2^*\rangle)\in
L^1(E,\mu_\infty)$, and $L^1(\R^2,\nu_\infty)$-convergence is
equivalent to $L^1(E,\mu_\infty)$-convergence for the corresponding
$(\langle x,h_1^*\rangle, \langle x,h_2^*\rangle)$-cylindrical
functions.

Now, for $\phi\in \mathcal{C}(L_2)$, we have $f(x):=\phi(\langle
x,h_1^*\rangle, \langle x,h_2^*\rangle)\in \mathcal{C}(L)$ and
\begin{eqnarray*}
\begin{aligned}
\langle &Ax,Df(x)\rangle=\langle Ax,\frac{\partial\phi}{\partial
t_1}\left(\langle x,h_1^*\rangle, \langle
x,h_2^*\rangle\right)\,h_1^*\rangle+\langle
Ax,\frac{\partial\phi}{\partial t_2}(\langle x,h_1^*\rangle, \langle
x,h_2^*\rangle)\,h_2^*\rangle\\&= \frac{\partial\phi}{\partial
t_1}(\langle x,h_1^*\rangle, \langle x,h_2^*\rangle)\,(a\langle
x,h_1^*\rangle-b\langle x,h_2^*\rangle)+
\frac{\partial\phi}{\partial t_2}(\langle x,h_1^*\rangle, \langle
x,h_2^*\rangle)\,(b\langle x,h_1^*\rangle+a\langle x,h_2^*\rangle),
\end{aligned}
\end{eqnarray*}
and, for $y\in H$,
\begin{eqnarray*}
\begin{aligned}
\left(D_H^2f(x)\right)&(y)=\frac{\partial^2 \phi}{\partial
t_1^2}(\langle x,h_1^*\rangle, \langle x,h_2^*\rangle) \langle
y,B^*h_1^*\rangle B^*h_1^*+\frac{\partial^2 \phi}{\partial t_1
\partial t_2}(\langle x,h_1^*\rangle, \langle x,h_2^*\rangle) \langle
y,B^*h_2^*\rangle B^*h_1^*\\& +\frac{\partial^2 \phi}{\partial t_2
\partial t_1}(\langle
x,h_1^*\rangle, \langle x,h_2^*\rangle)\langle y,B^*h_1^*\rangle
B^*h_2^* +\frac{\partial^2 \phi}{\partial t_2^2 }(\langle
x,h_1^*\rangle, \langle x,h_2^*\rangle) \langle y,B^*h_2^*\rangle
B^*h_2^*.
\end{aligned}
\end{eqnarray*}

Taking into account Lemma~\ref{trace}, we get
\begin{multline}\label{el2}
Lf(x)=\frac12\sum_{i,j=1}^2{r_{ij}\,\frac{\partial^2 \phi}{\partial
t_i
\partial t_j}(\langle
x,h_1^*\rangle, \langle x,h_2^*\rangle)}+ \sum_{i,j=1}^2{c_{ij}\,\,
\langle x,h_j^*\rangle\,\, \frac{\partial \phi}{\partial t_i
}(\langle x,h_1^*\rangle, \langle
x,h_2^*\rangle)}\\=(L_2\phi)(\langle x,h_1^*\rangle, \langle
x,h_2^*\rangle)).
\end{multline}

Now observe that $\sigma(C)=\{a\pm ib\}\subset \mathbb{C}_-$, and
the kernel of $R$ does not contain any invariant subspace of $C^*$:
by~(\ref{eq5}) this would imply degeneracy of $R_\infty$ and
$Q_\infty$, and consequently of $\mu_\infty$. Thus we can
use~\cite[Thm 5.1]{{M,P,Pr}} to conclude that for any
$\lambda\in\mathbb{C}_-$ there exists an eigenfunction
$\phi_\lambda(t)\in L^1(\R^2,\nu_\infty)$ of $L_2$ corresponding to
$\lambda$. Approximating $\phi_\lambda$ by
$\phi_n\in\mathcal{C}(L_2)$, $\phi_n\to\phi_\lambda$ in
$L^1(\R^2,\nu_\infty)$, and then applying~(\ref{el2}) to
$f_n(x):=\phi_n(\langle x,h_1^*\rangle, \langle
x,h_2^*\rangle)\in\mathcal{C}(L)$, we obtain $f_n(x)\to
\phi_\lambda(\langle x,h_1^*\rangle, \langle
x,h_2^*\rangle)=:f_\lambda(x)$ in $L^1(E,\mu_\infty)$, $Lf_n(x)\to
\lambda f_\lambda(x)$ in $L^1(E,\mu_\infty)$. Hence $f_\lambda(x)\in
\mathcal{D}(L)$ and $Lf_\lambda(x)=\lambda f_\lambda(x)$.
\end{proof}
\end{lemma}

\end{section}

\bigskip

\textsc{Dept of Mathematics, California Institute of Technology, MC
253-37, 1200 E.Ca\-li\-for\-nia Blvd, Pa\-sa\-de\-na, CA 91125, US.}

\textsc{Telephone: (626)29-89-119}

\textsc{E-mail: rostysla@caltech.edu}

\end{document}